\pdfoutput=1
\documentclass[a4paper]{article}
\usepackage{RR}
\usepackage{hyperref}
\usepackage{amsmath}
\usepackage{amssymb}
\usepackage{latexsym}
\RRdate{April 2008}
\RRNo{6490}
\RRauthor{% les auteurs
Ma\"elle Nodet\thanks{Universit\'e de Grenoble, INRIA. maelle.nodet@inria.fr}%
}
\authorhead{Nodet}
\RRtitle{Optimal control of the Primitive Equations of the Ocean with Lagrangian observations}
\RRetitle{Optimal control of the Primitive Equations of the Ocean with Lagrangian observations}
\titlehead{Optimal control of the Primitive Equations with Lagrangian observations}
\RRresume{On consid\`ere un probl\`eme de contr\^ole optimal pour les \'equations primitives non lin\'eaires de l'oc\'ean en dimension trois d'espace, dans un domaine verticalement born\'e et horizontalement p\'eriodique. L'op\'erateur d'observation associe \`a une solution des \'equations primitives la trajectoire d'une particule lagrangienne transport\'ee par le flot \`a profondeur donn\'ee. Ce travail montre l'existence d'un contr\^ole optimal, apr\`es r\'egularisation du probl\`eme. Pour cela, on \'etablit en particulier de nouvelles estimations d'\'energie pour les \'equations primitives dans des espaces fonctionnels bien choisis, qui distinguent la dimension verticale des dimensions horizontales. On illustre le r\'esultat par une exp\'erience num\'erique. 
}
\RRabstract{We consider an optimal control problem for the three-dimensional non-linear Primitive Equations of the ocean in a vertically bounded and horizontally periodic domain. The observation operator maps a solution of the Primitive Equations to the trajectory of a Lagrangian particle. This paper proves the existence of an optimal control for the regularized problem. To do that, we prove also new energy estimates for the Primitive Equations, thanks to well-chosen functional spaces, which distinguish the vertical dimension from the horizontal ones. We illustrate the result with a numerical experiment.}
\RRmotcle{contr\^ole optimal, \'equations aux d\'eriv\'ees partielles, \'equations primitives, simulation num\'erique}
\RRkeyword{optimal control, partial differential equations, Primitive Equations, numerical simulation}
\RRprojet{MOISE}  % cas d'un seul projet
\RRtheme{\THNum} % cas d'un seul theme
 \RCGrenoble % Grenoble - Rh\^one-Alpes
\newtheorem{prblm}{Problem}
\newtheorem{thrm}{Theorem}
\newtheorem{prpstn}{Proposition}

\newtheorem{rmrk}{Remark}
\newtheorem{dfntn}{Definition}

\newcommand{\bproof}{\paragraph{Proof. }}
\newcommand{\eproof}{$\Box$}

\begin{document}
\makeRR   % cas d'un rapport de recherche
%% \makeRT % cas d'un rapport technique.
%% a partir d'ici, chacun fait comme il le souhaite

\tableofcontents

\newpage

%
%%-----------------------------
\section*{Introduction}
%%-----------------------------
%

The ocean plays a major role in governing the earth climate. Physical oceanographers and climatologists work toward a better knowledge of the ocean properties (currents, temperature, salinity, marine biology, etc.). Some mathematical tools are involved in this progress, in particular Data Assimilation methods. Data Assimilation covers all mathematical methods which allow to blend optimally all sources of information about the ocean (measures, model equations, errors statistics) in order to improve ocean modeling, forecasts or climatology. One class of data assimilation methods, called Variational Data Assimilation \cite{LeDimet82,LeDimet86}, is based on optimal control theory \cite{Lions68}. The idea is to compute numerically the solution of an optimal control problem, in which the cost function represents the misfit between the observations and their model counterpart.\\
This paper deals with a particular problem of data assimilation for the ocean, namely the assimilation of Lagrangian data. The ocean is mainly observed at the surface, thanks to observing satellites. In-situ data are sparsely sampled in time and space, and it is therefore important to make the most of their information. Lagrangian data consist of positions of floats drifting at depth (around one thousand meters deep), they give information about in-depth currents. The problem of variational assimilation of Lagrangian data has been studied numerically in \cite{Nodet06}.\\
In this paper we investigate the theoretical justification of this problem. We prove the existence of an optimal control for Lagrangian observations for the Primitive Equations (PEs) of the ocean. To do so, we had to establish local existence and unicity of strong solutions for the PEs, allowing existence of Lagrangian trajectories, so that the problem can be formulated. The problem of local existence and unicity of strong solution for the PEs has been studied by Lions, Temam, Wang and Ziane \cite{LionsTemam92,TemamZiane04}. Due to the dissymmetry between the vertical dimension and the horizontal ones in our chosen domain (vertically bounded and horizontally periodic) and also in the PEs (as the equation for the vertical velocity is degenerated, contrary to Navier-Stokes equations), we had to introduce new functional spaces. In these spaces we successively prove new energy estimates for the linear PEs and existence of strong solutions for the non-linear PEs. \\
This paper is organized as follows. In section \ref{sec:main} we state the equations, the functional spaces, the cost function, the optimal control problem and the main results of the article. In section \ref{sec:proofPE} we prove new energy estimates for the linear PEs and local existence and unicity of the non-linear PEs. In section \ref{sec:proofCTRL} we prove the existence of an optimal control. Finally, in section \ref{sec:num} we present a numerical illustration of this problem.

%
%%-----------------------------
\section{Statement of the problem and main results}
\label{sec:main}
%%-----------------------------
%
\subsection{The Primitive Equations of the ocean}
\label{2sec:1}
We consider the Primitive Equations of the ocean in a three dimensional domain (see \cite{LionsTemam92,TemamZiane04}):
\begin{equation}
\left\{ \begin{array}
{ll} 
\label{2eq:1}
\partial_t u - \nu \Delta u +(U.\nabla_2)u +w\partial_z u- \alpha v +  \partial_x p   =  0  &\textrm{ in } \Omega \times (0,T)\\
\partial_t v - \nu \Delta v +(U.\nabla_2)v +w\partial_z v + \alpha u +  \partial_y p  =  0 &\\
 \partial_z p  - \beta \theta  = 0 &\\ \medbreak
\partial_{t} \theta -\nu \Delta \theta + (U.\nabla_2) \theta + w \partial_z \theta  +\gamma w =  0  & \textrm{ in } \Omega \times (0,T)\\\medbreak
w(x,y,z)=-\int_0^z \partial_x u(x,y,z') + \partial_y v(x,y,z') \, dz'& \textrm{ in } \Omega \times (0,T)\\ 
U(t=0) = U_0, \qquad \theta(t=0) = \theta_0& \textrm{ in } \Omega
\end{array} \right.
\end{equation}
with\\
\indent - $\Omega = \mathbb{T}^2 \times(0,a)$ the physical domain, with $\mathbb{T}^2=(\mathbb{R}/2\pi\mathbb{Z})^2$ the bidimensional torus, such that the domain is periodic in the horizontal directions $x$ and $y$, vertically bounded in $z$, with fixed depth;\\
\indent - $(0, T)$ the time interval;\\
\indent - $U=(u,v)$ the horizontal velocity vector, $w$ the vertical velocity, $\theta$ the temperature and $p$ the pressure;\\
\indent - $U_0=(u_0,v_0)$ and  $\theta_0$ the initial conditions;\\
\indent - $\nabla_2=(\partial_x ,\partial_{y} )$ the horizontal 2D gradient operator, $(\nabla_2 .)$ the horizontal divergence operator, $\nabla =(\partial_x ,\partial_{y} ,\partial_{z} )$ the 3D gradient operator, $\Delta = \partial_{xx} +\partial_{yy} +\partial_{zz} $ the 3D Laplacian;\\
\indent - $\alpha$, $\nu$, $\gamma$, $\beta$ physical constants.\\
The boundary conditions are the following:
\begin{equation}
\left\{ \begin{array}
{l} 
\label{2eq:2} 
\medbreak
u,v,\theta  \textrm{ are periodic in } x,y\\
\medbreak
u=0,v=0,\theta=0 \textrm{ on } \mathbb{T}^2\times \{z=0,z=a\}\times (0,T) \\
\int_{z=0}^a \partial_x u + \partial_y v \, dz= 0 \quad \textrm{ on } \mathbb{T}^2\times (0,T)
\end{array} \right.
\end{equation}
We denote by  $X(t)=(u(t),v(t),\theta(t))$ the state vector of our system, and by  $X_0=(u_0,v_0,\theta_0)$ the initial state, which will be our control.\\
We will focus on smooth solutions of the Primitives Equations (PEs), so that the cost function (involving Lagrangian trajectories) can be defined. To this end we first define the following functional spaces, for $m\in \mathbb{N}$: 
$$ %
\begin{array}{rcl}
%\label{2eq:3} 
\medbreak
 L^2_z H^{m}_{xy} &= & \big\{ u \in L^2(\Omega)\textrm{ periodic in }x,y,\, \partial^{\alpha}_{x,y} u \in L^2(\Omega), \forall \alpha \in \mathbb{N}^2, |\alpha| \leq m\big\}\medbreak\\
\mathcal{U}^{m+1}&=& \big\{ X=(u,v,\theta) \in (L^2_z H^{m}_{xy})^3, \textrm{ periodic in } x,y,\smallbreak\\
&& \quad X=0 \textrm{ on } \mathbb{T}^2\times \{z=0,z=a\},\smallbreak\\
&&\quad \int_{z=0}^a \partial_x u + \partial_y v \, dz= 0 \textrm{ on } \mathbb{T}^2 \smallbreak\\
&&\quad \nabla X=(\nabla u,\nabla v,\nabla \theta) \in ((L^2_z H^{m}_{xy})^3)^3 \, \big\} \medbreak\\
\mathcal{H}^{m+1}&=& \big\{ u \in L^2_z H^{m}_{xy},\textrm{ periodic in } x,y,\smallbreak\\
&& \quad u=0 \textrm{ on } \mathbb{T}^2\times \{z=0,z=a\}, \smallbreak\\
&&\quad \nabla u \in (L^2_z H^{m}_{xy})^3 \, \big\} 
\end{array} 
$$ %
associated to the following scalar product and norms:
$$ %
\begin{array}{l} 
%\label{2eq:4}
 (u_1,u_2)_{L^2_z H^m_{xy}} =  \displaystyle \sum_{|\alpha | \leq m} \displaystyle \int_{\Omega}  \partial^{\alpha}_{x,y} u_1 \partial^{\alpha}_{x,y} u_2 \,\,dx\,dy\,dz\\
 (\nabla u_1,\nabla u_2)_{(L^2_z H^m_{xy})^3} =  \displaystyle \sum_{| \alpha | \leq m} \displaystyle \int_{\Omega}  \partial^{\alpha}_{x,y} \nabla u_1 .\partial^{\alpha}_{x,y} \nabla u_2 \,\,dx\,dy\,dz\\
\begin{array}{ll}
 (X_1,X_2)_{\mathcal{U}^{m+1}} = &(u_1,u_2)_{L^2_z H^m_{xy}} + (\nabla u_1,\nabla u_2)_{(L^2_z H^m_{xy})^3} \\
\qquad &+  (v_1,v_2)_{L^2_z H^m_{xy}}+ (\nabla v_1,\nabla v_2)_{(L^2_z H^m_{xy})^3} \\
\qquad &+  K (\theta_1,\theta_2)_{L^2_z H^m_{xy}}+ K (\nabla \theta_1,\nabla \theta_2)_{(L^2_z H^m_{xy})^3} \medbreak
\end{array}\\
 \|u\|^2_{L^2_z H^m_{xy}} = (u,u)_{L^2_z H^m_{xy}}\\
 \|X\|^2_{\mathcal{U}^{m+1}} = (X,X)_{\mathcal{U}^{m+1}}
\end{array} 
$$ %
($K$ is a ``large'' constant which will be set later).\\
We define on $(\mathcal{H}^{m+1})^3$ the same scalar product and norm as on $\mathcal{U}^{m+1}$.
\begin{rmrk}
%\label{2rmq:2}
The functional spaces $\mathcal{U}^{m+1}$ and $(\mathcal{H}^{m+1})^3$ are not interpolation spaces. In the sequel, $m$ is a \emph{fixed} integer larger or equal to 2.
\end{rmrk}
In this framework, we have the
\begin{thrm} 
\label{2prop:1}
Let $m\geq 2$ be an integer and $X_0=(u_0,v_0,\theta_0) \in \mathcal{U}^{m+1}$. If $K$ is large enough, there exists  $t^* >0$ with $t^*=t^*(\alpha,\beta,\gamma,\nu,\|X_0\|_{\mathcal{U}^{m+1}})$ and there exists a unique solution  $X(t)=(u(t),v(t),\theta(t))$ of the PEs (\ref{2eq:1}) with boundary conditions (\ref{2eq:2}) such that
$$ %
%\label{2eq:30}
X \in \mathcal{C}([0,t^*];\mathcal{U}^{m+1}),\quad \partial_t X \in L^2(0,t^*;L^2_z H^m_{xy})
$$ %
Moreover, we have:
\begin{equation}
\label{2eq:56}
\| X(t)\|_{\mathcal{U}^{m+1}}^2  +\frac{1}{\nu} \int_0^t \| \partial_t X(s)\|_{2,m}^2 \, ds  \quad \leq \quad \frac{M}{\delta} \|X_0\|_{\mathcal{U}^{m+1}}^2 
\end{equation}
for all $t\in[0,t^*]$, where $\delta$ depends on $t^*$.
\end{thrm}
This result is proven in section  \ref{sec:proofPE}.
\subsection{The Lagrangian observations and the cost function}
We can assume, without loss of generality, that there is only one drifting float, and that its position is observed at only one given time $t_1$. Its position 
$\xi (t)=(\xi^1(t),\xi^2(t))$ in the plane $z=z_0$ is solution of the following differential equation: 
\begin{equation}
\label{2eq:5}
\left\{ \begin{array}{rcl}
 \dfrac{d\xi }{dt}& =& U(t,\xi^1(t),\xi^2(t),z_0)\medbreak\\
 \xi (0)& =& \xi_0 \end{array} \right.
\end{equation}
The following proposition is an easy consequence of theorem \ref{2prop:1}:
\begin{prpstn} %
%\label{2prop:2}
Under the hypothesis of proposition  \ref{2prop:1}, the unique solution  $X$ of the PEs (\ref{2eq:1}) and (\ref{2eq:2}) is continuous in time and in $z$, Lipschitz in $(x,y)$. Moreover, for all   $\xi_0 \in \mathcal{T}^2$ and $z_0 \in [0,a]$, there exists a unique Lagrangian  trajectory, solution of equation (\ref{2eq:5}), associated to  $X$, $\xi_0$ and $z_0$.
\end{prpstn}
We then define the following cost function:
\begin{equation}
\begin{array}{ccccc}
\label{2eq:6}
\mathcal{J}(X_0) & = & \frac{1}{2} \|\xi (t_1) -d\|^2 &+ & \frac{\omega}{2} \| X_0 \|_{\mathcal{U}^{m+1}}^2\\
& =& \mathcal{J}^o(X_0)  &+&\omega \, \mathcal{J}^b(X_0)  
\end{array} 
\end{equation}
with:\\
\indent - $d=(d_1,d_2)$ the observation;\\
\indent - $m$ an integer, $m\geq 2$, $\omega$ a positive constant;\\\indent - $\|.\|$ the Euclidian norm in the 2D plane $z=z_0$.\\
The \emph{observation operator} is thus defined as follows:
\begin{equation}
\label{2eq:121}
\mathcal{G}(t_1; X_0)= \xi (t_1)
\end{equation}
where $\xi$ is defined by equation (\ref{2eq:5}) where the velocity field  $U=(u,v)$ is solution of the PEs (\ref{2eq:1})  and (\ref{2eq:2}) initialized with $X_0$.
\begin{rmrk} 
%\label{2rmq:1}
Contrary to the classical theory of J.-L. Lions \cite{Lions68}, the observation operator is non linear; moreover it is defined as a function of either the initial state $X_0$, or as a function of the complete velocity field  $\{U(t), t\in[0,t_1]\}$, and not only  $U(t_1)$.
\end{rmrk}
\subsection{Statement of the problem and main result}
The optimal control problem associated to the observation of Lagrangian data is the following:
\begin{prblm}
\label{2pb:1}
Let $d\in\mathbb{R}^2$ be an observation. We look for an optimal control  $X_0^* \in \mathcal{U}^{m+1}$ solution of the following minimization problem:
$$ %
%\label{2eq:29}
\mathcal{J}(X_0^*)= \inf_{X_0 \in \mathcal{U}^{m+1}} \mathcal{J}(X_0)
$$ %
where the cost function is defined by  (\ref{2eq:6}), the state equation by  (\ref{2eq:1},\ref{2eq:2}) and the observations by (\ref{2eq:5}).
\end{prblm}
The main result of this paper is the
\begin{thrm} 
\label{2thm:1}
There exists an optimal control $X_0^* \in \mathcal{U}^{m+1}$ solution of problem \ref{2pb:1}.
\end{thrm}
This result is proven in section \ref{sec:proofCTRL}.
%
%%-----------------------------
\section{Existence of strong solutions for the Primitive Equations of the ocean}  
\label{sec:proofPE}
%%-----------------------------
%
In this section we prove theorem \ref{2prop:1} in three steps: first we prove energy estimates for the linear PEs, then we prove estimates for the non linear terms of (\ref{2eq:1})  and then we prove theorem \ref{2prop:1}.
\subsection{Energy estimates for the linear Primitive Equations}
We consider the following linear Primitive Equations: 
\begin{equation}
\left\{ \begin{array}
{ll} 
\label{2eq:7}
\partial_t u - \nu \Delta u - \alpha v +  \partial_x p   =  F_1  &\textrm{ in } \Omega \times (0,T)\\
\partial_t v - \nu \Delta v + \alpha u +  \partial_y p  =  F_2 &\\
 \partial_z p  - \beta \theta  = 0 &\\ \medbreak
\partial_{t} \theta -\nu \Delta \theta +\gamma w =  F_3  & \textrm{ in } \Omega \times (0,T)\\\medbreak
w(x,y,z)=-\int_0^z \partial_x u(x,y,z') + \partial_y v(x,y,z') \, dz'& \textrm{ in } \Omega \times (0,T)\\ 
U(t=0) = U_0, \qquad \theta(t=0) = \theta_0& \textrm{ in } \Omega
\end{array} \right.  
\end{equation}
with the same notations as in section \ref{2sec:1} and boundary conditions (\ref{2eq:2}).\\
In the sequel we will use the following notations:
$$\begin{array}{rcl}
\iint f& :=&\int_0^T \int_{\Omega} f(t,x,y,z) \,dx\,dy\,dz \,dt\smallskip\\
\|f\|&:=&\|f\|_{L^2(\Omega)} \smallskip\\
\|f\|_{2,m} &:=& \|f\|_{L^2_z H^m_{xy}} \smallskip\\
\|f\|_{m} &:=& \|f\|_{H^m_{xy}} \smallskip\\
\|(f_1,f_2)\|_{2,m} &:= &\|(f_1,f_2)\|_{(L^2_z H^m_{xy})^2} \smallskip\\
\|(f_1,f_2,f_3)\|_{2,m} &:=& \|(f_1,f_2,f_3)\|_{(L^2_z H^m_{xy})^3}
\end{array}
$$
The following proposition holds true:
\begin{prpstn} 
\label{2lem:1}
For all $K$ large enough, for all $T>0$, there exist constants  $C_1(a,\nu,K,\gamma,\beta)$, $C_2(K,\nu)$, $C_3(a,\nu)$ and $C_4(\nu)$ such that, for all  $X_0 \in \mathcal{U}^{m+1}$, $F \in L^2(0,T;L^2_z H^m_{xy})$, the unique solution  $X(t)$ of the linear PEs (\ref{2eq:7}) satisfies:
$$ %
%\label{2eq:111}
 X(t) \in \mathcal{C}([0,T],\mathcal{U}^{m+1})
$$ %
Moreover, the following inequality holds true:
\begin{equation}
\begin{array}{l}
\label{2eq:28}
\| X(t)\|_{2,m}^2 + \| \nabla X(t)\|_{2,m}^2  +\frac{1}{\nu} \int_0^t \| \partial_t X(s)\|_{2,m}^2 \, ds  \\ 
\quad \quad \leq e^{C_1 t}\big(
C_2 \| X_0\|_{2,m}^2 +C_3 \| \nabla X_0\|_{2,m}^2
 + C_4 \int_0^{T}\|F(s)\|_{2,m}^2 \, ds \big)
\end{array}
\end{equation}
for all  $t \in [0,T]$.
\end{prpstn}
\bproof
Classical variational methods (see for example \cite{LionsTemam92,TemamZiane04}) prove that for $X_0$ given in  $\mathcal{U}^{m+1}$ and $F \in L^2(0,T;L^2_z H^m_{xy})$, there exists $X(t)$ at least in $L^2(0,T;\mathcal{V})\cap \mathcal{C}([0,T];\mathcal{H})$, where $\mathcal{V}$ and $\mathcal{H}$ are classical spaces (see \cite{Girault79,Temam84,TemamHu02}) defined as follows:
\begin{dfntn}
%\label{dfn:2} 
Let
$$ %
%\label{7eq:21}
\begin{array}{rcl}
E_1 &=& \{ U=(u,v)\in\mathcal{C}^\infty(\Omega)^2, u,v\textrm{ periodic in }x,y, \\
&&\quad u=0,v=0 \textrm{ on } \mathbb{T}^2\times \{z=0,z=a\}\\
&&\quad\int_{0}^a \partial_x u(x,y,z') +\partial_y v(x,y,z') \, dz' =0, \forall (x,y)\in \mathbb{T}^2 \}\medbreak\\
E_2 &=& \{ \theta\in\mathcal{C}^\infty(\Omega),  \theta \textrm{ periodic in } x,y, \\
&&\quad\theta=0 \textrm{ on } \mathbb{T}^2\times \{z=0,z=a\}\}
\end{array}
$$ %
Then $\mathcal{H}_1$ (respectively $\mathcal{H}_2$) is defined to be the closure of $E_1$ in $L^2(\Omega)^2$ (resp. $L^2(\Omega)$),  and $\mathcal{V}_1$ (resp. $\mathcal{V}_2$) is the closure of $E_1$ (resp. $E_2$) in $H^1(\Omega)^2$ (resp. $H^1(\Omega)$), and finally $\mathcal{H}=\mathcal{H}_1\times \mathcal{H}_2$,  $\mathcal{V}=\mathcal{V}_1\times \mathcal{V}_2$.
\end{dfntn}
Thus it suffices to prove (\ref{2eq:28}). To this end, we successively state four energy estimates: first an estimate of $\|X(t)\|_{L^2_{xyz}}$ and similarly of $\|X(t)\|_{L^2_zH^m_{xy}}$, then an estimate of $\|\nabla X(t)\|_{L^2_{xyz}}$ and similarly of $\|\nabla X(t)\|_{L^2_zH^m_{xy}}$.\\
To obtain energy estimates for $\|X(t)\|_{L^2_{xyz}}$, we multiply equations (\ref{2eq:7}) by  $u$, $v$, $w$, $K \theta$ and we integrate in space an time. Thus we have:
\begin{equation}
\begin{array}{ll}
\label{2eq:8}
 T_1 + T_2 + T_3 + T_4 + T_5 = T_6&\medbreak\\
T_1 = \iint \partial_t u u + \partial_t v v + K \partial_t \theta \theta&
T_2 = \iint -\nu \Delta u u  -\nu \Delta v v  -\nu K \Delta \theta \theta \smallskip\\
T_3 = \iint -\alpha v u + \alpha u v&
T_4 = \iint -\beta \theta w + K \gamma w \theta \smallskip\\
T_5 = \iint \partial_x p u + \partial_y p v + \partial_z p w&
T_6 = \iint F_1 u + F_2 v + K F_3 \theta
\end{array}
\end{equation}
We integrate by parts using conditions (\ref{2eq:2}) and we get:
\begin{equation}
\begin{array}{lll}
\label{2eq:9}
T_1 &=& \frac{1}{2} \big( \|u(t)\|^2 +  \|v(t)\|^2 + K \|\theta(t)\|^2 - \|u_0\|^2 - \|v_0\|^2 - \|\theta_0\|^2   \big)\smallskip\\
&=&  \frac{1}{2} \big( \|X(t)\|^2 - \|X_0\|^2   \big)\medskip\\
T_2 &=& \nu \int_0^t \| \nabla u(s) \|^2 + \| \nabla v(s) \|^2 + K \| \nabla \theta(s) \|^2 \, ds \smallskip\\
&=& \nu \int_0^t \|\nabla X(s)\|^2 \medskip\\
T_3 &=& 0\medskip\\
T_5 &=& - \iint p (\partial_x u + \partial_y v + \partial_z w) +\int_0^t \int_{\mathbb{T}^2} p (w|_{z=1}- w|_{z=0}) \,dx\,dy \, dt\smallskip\\
&=&0
\end{array}
\end{equation}
Then (\ref{2eq:8}) and (\ref{2eq:9}) give:
\begin{equation}
\label{2eq:10}
\|X(t)\|^2 + 2\nu \int_0^t \|\nabla X(s)\|^2 =
\|X_0\|^2 + 2 \iint ( F_1 u + F_2 v + K F_3 \theta )
-2 (K\gamma-\beta) \iint w \theta 
\end{equation}
Then we give a bound of the right hand side of equality (\ref{2eq:10}). First we establish the following useful inequality for  $\|w\|$:
\begin{equation}
\label{2eq:31}
\begin{array}{rcl}
 \|w\|^2 &= &\int_{\Omega} |\int_0^z \partial_x u + \partial_y v \, dz'|^2 \, dx\,dy\,dz \smallskip\\
&\leq& \int_{\Omega} 2 z \int_0^z |\partial_x u|^2 + |\partial_y v|^2 \, dz' \, dx\,dy\,dz \smallskip\\
&\leq& a^2 (\|\partial_x u \|^2 + \|\partial_y v\|^2)
\end{array}
\end{equation}
Thanks to (\ref{2eq:31}) we have:
\begin{equation}
\label{2eq:11}
\begin{array}{rcl}
2 \iint ( F_1 u + F_2 v + K F_3 \theta ) &\leq & \int_0^t \|F_1(s)\|^2 +\|F_2(s)\|^2 +K \|F_3(s)\|^2 \smallskip\\
&&+ \|u(s)\|^2 + \|v(s)\|^2 +K \|\theta(s)\|^2 \, ds \smallskip\\
&=& \int_0^t \|F(s)\|^2 + \|X(s)\|^2 \, ds \medskip\\
-2 (K\gamma-\beta) \iint w \theta &\leq & |K\gamma-\beta| \int_0^t (\|w\|^2 + \|\theta\|^2) \smallskip\\
&\leq& |K\gamma-\beta| \int_0^t  a^2 \|\partial_x u \|^2 + a^2  \|\partial_y v\|^2 + \|\theta\|^2 \, ds
\end{array}
\end{equation}
Then (\ref{2eq:10}) and (\ref{2eq:11}) lead to:
$$ %
\begin{array}{rcl}
%\label{2eq:12}
\|X(t)\|^2 + 2\nu \int_0^t \|\nabla X(s)\|^2 &\leq&\|X_0\|^2 + \int_0^t \|F(s)\|^2 + \int_0^t\|X(s)\|^2 \, ds \smallskip\\
&& + \int_0^t \|\partial_z u \|^2 + \|\partial_z v \|^2 \, ds \smallskip\\
&&+ |K\gamma-\beta| \int_0^t  a^2  \|\partial_x u \|^2 + a^2  \|\partial_y v\|^2 + \|\theta\|^2 \, ds \medskip\\
&\leq& \|X_0\|^2  + \int_0^t \|F(s)\|^2 + \int_0^t \|U(s)\|^2 \smallskip\\
&&  + (K + |K\gamma-\beta|) \|\theta(s)\|^2 \, ds \smallskip\\
&& + \int_0^t \max(1, a^2 |K\gamma-\beta|)\|\nabla U(s)\|^2 \, ds
\end{array}
$$ %
We proceed similarly to obtain a bound on $\|X(t)\|_{L^2_zH^m_{xy}}$. First we note that if  $u$, $v$, $w$, $\theta$ and $p$ satisfy equation  (\ref{2eq:7}), then for all $\alpha \in \mathbb{N}^2$, $\partial_{xy}^{\alpha} u$, $\partial_{xy}^{\alpha} v$, $\partial_{xy}^{\alpha} w$, $\partial_{xy}^{\alpha} \theta$ and $\partial_{xy}^{\alpha} p$ satisfy the same equation  (where $F_i$ is replaced by $\partial_{xy}^{\alpha} F_i$ and $X_0$ by $\partial_{xy}^{\alpha} X_0$) and the same boundary conditions. Therefore we can apply the same calculations and obtain the following inequality:
\begin{equation}
\begin{array}{rcl}
\label{2eq:13}
\|X(t)\|_{2,m}^2 + 2\nu \int_0^t \|\nabla X(s)\|_{2,m}^2 & \leq &  \|X_0\|_{2,m}^2 + \int_0^t \|F(s)\|_{2,m}^2   \smallskip\\
&&  + \int_0^t \|U(s)\|_{2,m}^2 + (K + |K\gamma-\beta|) \|\theta(s)\|_{2,m}^2 \, ds \smallskip\\
&& + \int_0^t \max(1, a^2 |K\gamma-\beta|)\|\nabla U(s)\|_{2,m}^2 \, ds
\end{array}
\end{equation}
To establish the second type of estimates, we then multiply the equation by  $\partial_t u$, $\partial_t v$, $\partial_t w$ and $K \partial_t \theta$ and we integrate over $\Omega \times (0,t)$:
$$ %
\begin{array}{l}
%\label{2eq:14}
\medbreak
 T_1 + T_2 + T_3 + T_4 + T_5 = T_6\\
T_1 = \iint \partial_t u \partial_t u + \partial_t v \partial_t v + K \partial_t \theta \partial_t \theta \smallskip\\
T_2 = \iint -\nu \Delta u \partial_t u  -\nu \Delta v \partial_t v  -\nu K \Delta \theta \partial_t \theta \smallskip\\
T_3 = \iint -\alpha v \partial_t u + \alpha u \partial_t v \smallskip\\
T_4 = \iint -\beta \theta \partial_t w + K \gamma w \partial_t \theta \smallskip\\
T_5 = \iint \partial_x p \partial_t u + \partial_y p \partial_t v + \partial_z p \partial_t w \smallskip\\
T_6 = \iint F_1 \partial_t u + F_2 \partial_t v + K F_3 \partial_t \theta
\end{array}
$$ %
Using integration by parts and boundary conditions (\ref{2eq:2}) we get:
$$ %
\begin{array}{lll}
%\label{2eq:15}
T_1 &=& \int_0^t \|\partial_t u(s) \|^2  + \|\partial_t v(s) \|^2  +K \|\partial_t u(s) \|^2 \, ds \smallskip\\
    &=&  \int_0^t \| \partial_t X(s)\|^2 \, ds \medskip\\
T_2 &=& \frac{\nu}{2} \|\nabla X(t)\|^2 - \frac{\nu}{2} \|\nabla X_0\|^2 \medskip\\
T_4 &=& \iint (K\gamma +\beta) w \partial_t \theta + \int_{\Omega} \big( \theta_0 w_0 - \theta(t) w(t)\big) \, dx\,dy\,dz \medskip\\
T_5 &=& 0
\end{array}
$$ %
Thus:
$$ %
\begin{array}{rcl}
%\label{2eq:16}
2 \int_0^t \| \partial_t X(s)\|^2 \, ds
 + \nu \|\nabla X(t)\|^2 
&=&  \nu \|\nabla X_0\|^2 + 2 \iint \alpha v \partial_t u - \alpha u \partial_t v
\smallskip\\
&& 
- 2 \iint (K\gamma +\beta) w \partial_t \theta + 2 \int_{\Omega} (  \theta(t) w(t)- \theta_0 w_0) \, dx\,dy\,dz \smallskip\\
&&  + 2 \iint F_1 \partial_t u + F_2 \partial_t v + K F_3 \partial_t \theta
\end{array}
$$ %
In order to get a bound on the right hand side we use the following inequality:
$$ %
%\label{2eq:18}
xy \leq \frac{\varepsilon}{2}x^2 + \frac{1}{2\varepsilon}y^2 \qquad \forall x, y \in \mathbb{R}, \forall \varepsilon > 0
$$ %
Thus we obtain, for all positive real numbers  $\varepsilon_i$, $2\leq i \leq 5$:
$$ %
\begin{array}{rcl}
%\label{2eq:17}
2 \iint \alpha v \partial_t u - \alpha u \partial_t v &\leq & \alpha \int_0^t ( \varepsilon_2 \|\partial_t U(s)\|^2 + \frac{1}{\varepsilon_2} \|U(s)\|^2)\, ds \medskip\\
- 2 \iint (K\gamma +\beta) w \partial_t \theta &\leq &(K\gamma +\beta) \int_0^t \varepsilon_3 \|\partial_t \theta(s)\|^2 \, ds + \smallskip\\
&&(K\gamma +\beta) \int_0^t \frac{ a^2 }{\varepsilon_3} (\|\partial_x u(s) \|^2+\|\partial_y v(s) \|^2)\, ds \medskip\\
2 \int_{\Omega} (  \theta(t) w(t)- \theta_0 w_0) \, dx\,dy\,dz &\leq & a^2 \varepsilon_4 (\|\partial_x u(t)\|^2 + \|\partial_y v(t)\|^2) + \frac{1}{\varepsilon_4} \|\theta(t)\|^2 \smallskip\\
&&+\|\theta_0\|^2 +a^2 \|\partial_x u_0\|^2 + a^2 \|\partial_y v_0\|^2 \medskip\\
2 \iint F_1 \partial_t u + F_2 \partial_t v + K F_3 \partial_t \theta &\leq& \int_0^t(\frac{1}{\varepsilon_5}\|F(s)\|^2 + \varepsilon_5 \|\partial_t X(s)\|^2)\, ds
\end{array}
$$ %
And finally we get:
$$ %
\begin{array}{rcl}
%\label{2eq:19}
(\nu-a^2\varepsilon_4)\|\nabla U(t)\|^2 + K\nu \|\nabla \theta(t)\|^2 && \smallskip\\
+ \int_0^t (2-\alpha\varepsilon_2-\varepsilon_5)\| \partial_t U(s)\|^2  - \frac{1}{\varepsilon_4} \|\theta(t)\|^2 && \smallskip\\
 + (2K-(K\gamma+\beta)\varepsilon_3-K\varepsilon_5)\| \partial_t \theta(s)\|^2\, ds&
\leq &  (\nu+a^2) \|\nabla U_0\|^2 + K\nu\|\nabla \theta_0\|^2
   +\|\theta_0\|^2 \smallskip\\ 
&&+(K\gamma +\beta) \int_0^t \frac{ a^2 }{\varepsilon_3} \|\nabla U(s) \|^2\, ds+
\alpha \int_0^t   \frac{1}{\varepsilon_2} \|U(s)\|^2 \, ds \smallskip\\
&&+ 
\int_0^t\frac{1}{\varepsilon_5} \|F(s)\|^2 \, ds\end{array}
$$ %
As previously we obtain the same result for the derivative in $x$ and $y$:
$$ %
\begin{array}{rcl}
%\label{2eq:20}
(\nu-a^2\varepsilon_4)\|\nabla U(t)\|_{2,m}^2 + K\nu \|\nabla \theta(t)\|_{2,m}^2 && \smallskip\\
 - \frac{1}{\varepsilon_4} \|\theta(t)\|_{2,m}^2 
+ \int_0^t (2-\alpha\varepsilon_2-\varepsilon_5)\| \partial_t U(s)\|_{2,m}^2 && \smallskip\\
 + (2K-(K\gamma+\beta)\varepsilon_3-K\varepsilon_5)\| \partial_t \theta(s)\|_{2,m}^2\, ds&
 \leq  & (\nu+a^2) \|\nabla U_0\|_{2,m}^2 + K\nu\|\nabla \theta_0\|_{2,m}^2
   +\|\theta_0\|_{2,m}^2 \smallskip\\ 
&& +(K\gamma +\beta) \int_0^t \frac{ a^2 }{\varepsilon_3} \|\nabla U(s) \|_{2,m}^2\, ds \smallskip\\
&&+
\alpha \int_0^t   \frac{1}{\varepsilon_2} \|U(s)\|_{2,m}^2 \, ds+ 
\int_0^t\frac{1}{\varepsilon_5} \|F(s)\|_{2,m}^2 \, ds
\end{array}
$$ %
The $\varepsilon_i$ are chosen as follows:
$$ %
%\label{2eq:21}
\varepsilon_2 = \frac{1}{\alpha}, \quad
\varepsilon_3 = \frac{1}{\gamma}, \quad
\varepsilon_4 = \frac{\nu}{2a^2}, \quad
\varepsilon_5 = \frac{1}{2}
$$ %
and thus we get:
\begin{equation}
\begin{array}{rcl}
\label{2eq:22}
\frac{\nu}{2}\|\nabla U(t)\|_{2,m}^2 + K\nu \|\nabla \theta(t)\|_{2,m}^2 && \smallskip\\
 - \frac{2a^2}{\nu} \|\theta(t)\|_{2,m}^2 + \int_0^t \frac{1}{2}\| \partial_t U(s)\|_{2,m}^2 && \smallskip\\
 + (\frac{K}{2}-\gamma\beta)\| \partial_t \theta(s)\|_{2,m}^2\, ds& \leq  & (\nu+a^2) \|\nabla U_0\|_{2,m}^2 + K\nu\|\nabla \theta_0\|_{2,m}^2
  +\| \theta_0\|_{2,m}^2\smallskip\\ 
&& + \int_0^t  a^2 \gamma(K\gamma +\beta) \|\nabla U(s) \|_{2,m}^2\, ds \smallskip\\
&&+
 \int_0^t  \alpha^2 \|U(s)\|_{2,m}^2 \, ds
 + 
\int_0^t 2 \|F(s)\|_{2,m}^2 \, ds\end{array}
\end{equation}
We then add equation (\ref{2eq:13}) and equation (\ref{2eq:22}) multiplied by $\frac{2}{\nu}$ to obtain:
$$ %
\begin{array}{l}
%\label{2eq:23}
\| U(t)\|_{2,m}^2 +   (K-\frac{4a^2}{\nu^2}) \|\theta(t)\|_{2,m}^2
+ \|\nabla U(t)\|_{2,m}^2 +  2K \|\nabla \theta(t)\|_{2,m}^2   \smallskip\\
+ \int_0^t \frac{1}{\nu}\| \partial_t U(s)\|_{2,m}^2 
+ (\frac{K}{2}-\gamma\beta)\frac{2}{\nu}\| \partial_t \theta(s)\|_{2,m}^2\, ds
+ \int_0^t 2\nu\| \nabla U(s)\|_{2,m}^2 
+ 2\nu K\| \nabla \theta(s)\|_{2,m}^2\, ds\medskip \\
\qquad \qquad \leq   \qquad \| U_0\|_{2,m}^2 + (K+\frac{2}{\nu})\| \theta_0\|_{2,m}^2
 +  (2+\frac{2a^2}{\nu})\|\nabla U_0\|_{2,m}^2 +2K \|\nabla \theta_0\|_{2,m}^2 \smallskip\\
\qquad \qquad \qquad   \quad  + \int_0^t  (1+\frac{4}{\nu})\|F(s)\|_{2,m}^2 \, ds 
  \int_0^t  (1+\frac{2\alpha^2}{\nu}) \|U(s)\|_{2,m}^2 \, ds  
  + \int_0^t  (K+|K\gamma-\beta|) \|\theta(s)\|_{2,m}^2 \, ds \smallskip\\
\qquad \qquad \qquad   \quad + \int_0^t  (\max(1,a^2 |K\gamma-\beta|)+\frac{2\gamma a^2 }{\nu}(K\gamma+\beta))\|\nabla U(s) \|_{2,m}^2\, ds
\end{array}
$$ %
For $K \geq 2\max(\frac{4a^2}{\nu^2},2\gamma\beta)$, for all $T$ and for all $t\in[0,T]$ we have:
$$ %
%\label{2eq:24}
\| X(t)\|_{2,m}^2 + \| \nabla X(t)\|_{2,m}^2 +\frac{1}{\nu} \int_0^t \| \partial_t X(s)\|_{2,m}^2 \, ds
 \leq C_5 + C_1 \int_0^t \| X(s)\|_{2,m}^2 + \| \nabla X(s)\|_{2,m}^2\, ds
$$ %
with:
$$ %
\begin{array}{rl}
%\label{2eq:25}
 C_5 = &2\| U_0\|_{2,m}^2 + (2K+\frac{4}{\nu})\| \theta_0\|_{2,m}^2
+  (4+\frac{4a^2}{\nu})\|\nabla U_0\|_{2,m}^2 +4K \|\nabla \theta_0\|_{2,m}^2 \smallskip\\
& + 
\int_0^{T}  (2+\frac{8}{\nu})\|F(s)\|_{2,m}^2 \, ds\medbreak\\
 C_1 = & 2 \max\big(1+\frac{2\alpha^2}{\nu},\frac{1+|\gamma-\beta/K|}{1+ 2\alpha^2/\nu},\max(1,a^2 |K\gamma-\beta|)+\frac{2\gamma a^2 }{\nu}(K\gamma+\beta)\big)
\end{array}
$$ %
With Gronwall lemma we get:
$$ %
%\label{2eq:26}
C_5 + C_1 \int_0^t \| X(s)\|_{2,m}^2 + \| \nabla X(s)\|_{2,m}^2\, ds \leq C_5 e^{C_1 t}
$$ %
thus
$$ %
\begin{array}{l}
%\label{2eq:27}
\| X(t)\|_{2,m}^2 + \| \nabla X(t)\|_{2,m}^2  +\frac{1}{\nu} \int_0^t \| \partial_t X(s)\|_{2,m}^2 \, ds  \smallskip\\ 
\qquad \qquad\qquad \qquad\qquad \qquad \leq e^{C_1 t}\big(
C_2 \| X_0\|_{2,m}^2 +C_3 \| \nabla X_0\|_{2,m}^2
 + C_4 \int_0^{T}\|F(s)\|_{2,m}^2 \, ds \big)
\end{array}
$$ %
with
$$C_2 = 2+\frac{4}{K\nu}, \,\, C_3=4+\frac{4a^2}{\nu} \textrm{ et } C_4 = 2+\frac{8}{\nu}$$
\eproof
\subsection{Estimation of the non linear terms}
We will now estimate the non linear terms of equation  (\ref{2eq:1}). The following proposition holds true:
\begin{prpstn} 
\label{2lem:2}
Let $m\geq2$ be an integer, $X_1$ and $X_2$ elements of  $\mathcal{U}^{m+1}$. Let us define: 
$$ %
\begin{array}{rcl}
%\label{2eq:32}
F_1 &=& (U_1 . \nabla_2) u_2 + w_1 \partial_z u_2\\
F_2 &=& (U_1 . \nabla_2) v_2 + w_1 \partial_z v_2\\
F_3 &=& (U_1 . \nabla_2) \theta_2 + w_1 \partial_z \theta_2\\
\end{array}
$$ %
then we have, for all  $i\in\{1,2,3\}$:
$$ %
%\label{2eq:44}
\|F_i\|^2_{2,m} \quad \leq \quad C \,\,(\|X_1\|_{2,m}+a^2\|\nabla X_1 \|_{2,m}) \,\,\|\nabla X_1 \|_{2,m}\,\,\|\nabla X_2 \|^2_{2,m} 
$$ %
where $C$ is a constant (independent of $a$) of order $1$.\\
If $X_1=X_2$ then we have, for all $i\in\{1,2,3\}$:
$$ %
%\label{2eq:45}
\|F_i\|^2_{2,m}\quad \leq \quad C_5 \,\, \big( \|X_1\|^2_{2,m}+\|\nabla X_1 \|^2_{2,m} \big)^2
$$ %
where $C_5 = C_5(a)$ is a constant such that  $C_{5} = C'_5 + C''_5 a^2$ (where $C'_5$ and $C''_5$ are constants of order $1$).
\end{prpstn}
\bproof
First let us address the term $u_1\partial_x u_2$:
\begin{equation}
\begin{array}{rcl}
\label{2eq:33}
\| u_1\partial_x u_2 \|^2_{2,m} &=& \int_{z=0}^a \|u_1(z)\partial_x u_2(z) \|^2_{m}\,dz\smallbreak\\
&\leq& \int_{z=0}^a \|u_1(z) \|^2_{m}\|\partial_x u_2(z) \|^2_{m}\,dz
\end{array}
\end{equation}
because $H^m_{xy}$ is an algebra for  $m\geq 2$.\\
We estimate now $\sup_{z\in(0,a)} \|u_1(z) \|^2_{m}$. We first write an estimation in space dimension 1 for  $\varphi(z) \in H^1(0,a)$ with $\varphi(0)=0$:
$$ %
%\label{2eq:34}
\begin{array}{rcl}
|\varphi(z)|^2 &=& \int_0^z \partial_z |\varphi(z')|^2\,dz'\\
&=&\int_0^z 2\varphi(z')\partial_z \overline{\varphi}(z')\,dz' \\
&\leq &2 \big( \int_0^a |\varphi(z)|^2\,dz \,\int_0^a |\partial_z \varphi(z)|^2\,dz \big)^{1/2}
\end{array}
$$ %
Similarly we can write an estimation for $\|u_1(z) \|^2_{m}$:
\begin{equation}
\label{2eq:35}
\| u_1(z) \|^2_{m} \leq 2 \| u_1 \|_{2,m} \|\partial_z u_1 \|_{2,m}
\end{equation}
Therefore (\ref{2eq:33}) and (\ref{2eq:35}) give us:
$$ %
%\label{2eq:36}
\begin{array}{rcl}
\| u_1\partial_x u_2 \|^2_{2,m} &\leq&
2 \| u_1 \|_{2,m} \|\partial_z u_1 \|_{2,m} \int_{z=0}^a \|\partial_x u_2(z) \|^2_{m}\,dz\\
&\leq& 2 \|X_1\|_{2,m}\|\nabla X_1 \|_{2,m} \|\nabla X_2 \|^2_{2,m} 
\end{array}
$$ %
And if $u_1=u_2$ we have:
$$ %
%\label{2eq:37}
\begin{array}{rcl}
\| u_1\partial_x u_1 \|^2_{2,m} 
&\leq& 2 \|X_1\|_{2,m} \|\nabla X_1 \|^3_{2,m} \\
&\leq& \big( \|X_1\|^2_{2,m}+\|\nabla X_1 \|^2_{2,m} \big)^2
\end{array}
$$ %
We establish the same inequalities for the terms  $u_1 \partial_x u_2$, $u_1 \partial_x \theta_2$, $v_1\partial_y u_1$, $v_1\partial_y v_2$ and $v_1\partial_y \theta_2$.\\
Let us now focus on  $w_1 \partial_z u_2$. As in (\ref{2eq:33}), we have:
$$ %
%\label{2eq:38}
\| w_1\partial_z u_2 \|^2_{2,m} \leq \int_{z=0}^a \|w_1(z) \|^2_{m}\|\partial_z u_2(z) \|^2_{m}\,dz
$$ %
We have also, as in (\ref{2eq:35}):
$$ %
%\label{2eq:39}
\| w_1(z) \|^2_{m} \leq 2 \| w_1 \|_{2,m} \|\partial_z w_1 \|_{2,m}
$$ %
Thus:
$$ %
%\label{2eq:40}
\begin{array}{rcl}
\| w_1\partial_z u_2 \|^2_{2,m} &\leq& 2 \| w_1 \|_{2,m} \|\partial_z w_1 \|_{2,m}\|\partial_z u_2\|^2_{2,m}\smallbreak\\
&\leq& C \| w_1 \|_{2,m} \|\nabla X_1 \|_{2,m}\|\nabla X_2\|^2_{2,m}
\end{array}
$$ %
Then we can bound  $\| w_1 \|_{2,m}$ as in (\ref{2eq:31}):
$$ %
%\label{2eq:41}
\begin{array}{rcl}
\| w_1 \|^2_{2,m} &=& \sum_{|\alpha|\leq m }\int_0^a | \partial_{xy}^{\alpha} w_1|^2\, dz\smallbreak\\
 &=& \sum_{|\alpha|\leq m }\int_0^a | \partial_{xy}^{\alpha} \int_0^z \partial_x u_1(z')+\partial_y v_1(z') \,dz'|^2\, dz\smallbreak\\
&\leq& \sum_{|\alpha|\leq m }\int_0^a  z\, dz \int_0^a 2 | \partial_{xy}^{\alpha}  \partial_x u_1(z)|^2 +2 | \partial_{xy}^{\alpha} \partial_y v_1(z)|^2 \, dz\smallbreak\\
&\leq& a^2 \big( \|\partial_x u_1 \|_{2,m}^2 + \|\partial_y v_1\|_{2,m}^2 \big)\smallbreak\\
&\leq& a^2 \|\nabla X_1\|_{2,m}^2
\end{array}
$$ %
therefore we get:
$$ %
%\label{2eq:42}
\| w_1\partial_z u_2 \|^2_{2,m} \quad \leq\quad  C a^2 \|\nabla X_1 \|^2_{2,m}\|\nabla X_2\|^2_{2,m}
$$ %
And if $u1=u_2$ we have: 
$$ %
%\label{2eq:43}
\| w_1\partial_z u_1 \|^2_{2,m} \quad \leq\quad  C a^2 \|\nabla X_1 \|^4_{2,m} \quad \leq\quad  C a^2 \big( \|X_1\|^2_{2,m}+\|\nabla X_1 \|^2_{2,m} \big)^2
$$ %
We have the same estimates for $w_1 \partial_z v_2$ et $w_1 \partial_z \theta_2$, and that concludes the proof of proposition \ref{2lem:2}.\\
\eproof
\subsection{End of the proof}
We will now construct a solution of the non linear equation (\ref{2eq:1}). Let us first introduce some notations. \\
Let $L$ be an operator defined as follows:
$$ %
%\label{2eq:112}
L(X,p_0) = \left[ \begin{array}{l}
\partial_t u -\Delta u -\alpha v +\partial_x p_0 +\int_0^z \beta \partial_x \theta\,dz\smallbreak\\
\partial_t v -\Delta v +\alpha u +\partial_y p_0 +\int_0^z \beta \partial_y \theta\,dz\smallbreak\\
\partial_t \theta -\Delta \theta + \gamma w
\end{array} \right]
$$ %
And let us define $F(X_1,X_2)$ by:
$$ %
%\label{2eq:114}
F(X_1,X_2) = \left[ \begin{array}{l}
-(U_1.\nabla_2)u_2 - w_1 \partial_z u_2\\
-(U_1.\nabla_2)v_2 - w_1 \partial_z v_2\\
-(U_1.\nabla_2)\theta_2 - w_1 \partial_z \theta_2
\end{array} \right]
$$ %
We now define $N$, which is the square of a norm:
$$ %
%\label{2eq:115}
N(X(t)) =  \|X(t)\|^2 _{\mathcal{U}^{m+1}} + \int_0^t \|\partial_t X(s)\|^2_{2,m} \, ds 
$$ %
We are then given an initial condition $X_0\in\mathcal{U}^{m+1}$. Let us define the sequence $(X^n,p_0^n)$ in the following way:
$$ %
%\label{2eq:113}
\begin{array}{l}
X^0(t,x,y,z)  = X_0(x,y,z),\quad\forall (t,x,y,z)\in\mathbb{R}_+\times\Omega\medbreak\\
\left\{ \begin{array}
{l}
L(X^{n+1},p_0^{n+1}) = F(X^n,X^n)\qquad \textrm{pour }n\geq 0\smallbreak\\
 X^{n+1}|_{t=0} = X_0
\end{array} \right.
\end{array}
$$ %
We now verify that the sequence  $(X^n,p_0^n)$ is well defined for all $t$ and that $N(X^n(t))$ is finite for all $n$ and all $t$: we denote by  $N^0$ the value $N(X^0) = \|X_0\|^2 _{\mathcal{U}^{m+1}}$. Let us assume that  $X^n$ is well defined with a finite norm $N^{1/2}$  for all $t$. Thanks to propositions \ref{2lem:1} and \ref{2lem:2}, we get existence and unicity of  $X^{n+1}$ and moreover:
$$ %
%\label{2eq:116}
\begin{array}{rcl}
N(X^{n+1}(t))&\leq&C_0 e^{C_1 t}( N^0 + \int_0^{t} \| F(X^n,X^n)\|^2_{2,m} \, ds )\smallbreak\\
&\leq&C_0 e^{C_1 t} ( N^0 + C_5 \int_0^{t} \| X^n\|^4_{\mathcal{U}^{m+1}} \, ds )\smallbreak\\
&\leq&C_0 e^{C_1 t} ( N^0 + C_5 t \sup_{s\in[0,t]}N(X^n(s))^2  )\smallbreak\\
&\leq&C_6 e^{C_1 t} ( N^0 + t  \sup_{s\in[0,t]}N(X^n(s))^2 )\smallbreak\\
&<& \infty
\end{array}
$$ %
Now, let us choose $t^*\leq (4 C_6^2 e^{2 C_1 t^*} N^0)^{-1}$, we then have by recurrence 
$$ %
%\label{2eq:117}
\sup_{t\in[0,t^*]} N(X^{n}(t)) \leq 2 C_6 e^{ C_1 t^*} N^0
$$ %
In fact, we have $N_0 \leq 2 C_6 e^{ C_1 t^*} N^0 $ (after increasing $C_6$ if necessary) and then by recurrence:
$$ %
%\label{2eq:118}
\begin{array}{rcl}
\sup_{t\in[0,t^*]} N(X^{n+1}(t))&\leq&C_6 e^{C_1 t^*} ( N^0 + t^*  \sup_{s\in[0,t]}N(X^n(s))^2 )\smallbreak\\
&\leq&C_6 e^{C_1 t^*} ( N^0 + t^* (2 C_6 e^{ C_1 t^*} N^0)^2   )\smallbreak\\
&\leq&  2 C_6 e^{ C_1 t^*} N^0
\end{array}
$$ %
The sequence $(X^n)$ is therefore bounded for  $N$ uniformly in $t\in[0,t^*]$. \\
To pass to the limit, we write the equation verified by  $X^{n+1}-X^n$:
$$ %
%\label{2eq:119}
L(X^{n+1}-X^n) = F(X^n,X^n) -F(X^{n-1},X^{n-1}),\quad X^{n+1}-X^n|_{t=0} = 0
$$ %
Thanks to propositions \ref{2lem:1} and \ref{2lem:2} we obtain:
$$ %
%\label{2eq:120}
\begin{array}{rcl}
\sup_{t\in[0,t^*]} N(X^{n+1}-X^n)&\leq& C \int_0^{t^*} \| F(X^n,X^n) -F(X^{n-1},X^{n-1})\|^2_{2,m}\,dt\smallbreak\\
&=& C \int_0^{t^*} \|  \frac12 F(X^n-X^{n-1},X^n+X^{n-1}) \\
&& \qquad +  \frac12 F(X^n+X^{n-1},X^n-X^{n-1})  \|^2_{2,m}\,dt\smallbreak\\
&\leq& C \int_0^{t^*} N(X^n-X^{n-1})N(X^n+X^{n-1})\,dt\smallbreak\\
&\leq& C \int_0^{t^*} N(X^n-X^{n-1})(N(X^n)+N(X^{n-1}))\,dt\smallbreak\\
&\leq& C N^0 \int_0^{t^*} N(X^n-X^{n-1})\,dt\smallbreak\\
&\leq& C N^0 t^* \sup_{t\in[0,t^*]}N(X^n-X^{n-1})\smallbreak\\
&\leq& (C N^0 t^* )^n \sup_{t\in[0,t^*]}N(X^1-X^{0})
\end{array}
$$ %
We can again decrease  $t^*$ (if necessary) such that  $C N^0 t^*<1$, consequently   $(X^n)$ is a Cauchy sequence with respect to  $N$, and thus it converges to  $X$ strongly in $\mathcal{C}([0,t^*],\mathcal{U}^{m+1})$; $\partial_t X^n$ converges to  $\partial_t X$ strongly in  $L^2(0,t^*;L^2_zH^m_{xy})$ and inequality (\ref{2eq:56}) is true. \\
We can then pass to the limit into the variational formulation of the equation, and reintroduce the pressure, thanks to very classical methods (see for example \cite{LionsTemam92} and the further paper by Temam and Ziane \cite{TemamZiane04}), and this concludes the proof of  theorem \ref{2prop:1}.

%
%%-----------------------------
\section{Existence of an optimal control}
\label{sec:proofCTRL}
%%-----------------------------
%
In this section we prove theorem \ref{2thm:1}.\\
The proof uses the minimizing sequences method, in two steps: first we prove the convergence of the observation term, and then we pass to limit in the state equation.
\subsection{Convergence of the observation term}
Let $(X^n_0)$ be a minimizing sequence for $\mathcal{J}$:
$$ %
%\label{2eq:58}
\mathcal{J}(X^n_0) \rightarrow \inf_{X_0\in\mathcal{U}^{m+1}} \mathcal{J}(X_0) = \inf_{X_0\in\mathcal{U}^{m+1}} \mathcal{J}^o(X_0) + \|X_0\|^2_{\mathcal{U}^{m+1}}
$$ %
Then  $(X^n_0)$ is bounded in $\mathcal{U}^{m+1}$ and in $(\mathcal{H}^{m+1})^3$, and converges to  $X^*_0$ weakly in $\mathcal{U}^{m+1}$ and in $(\mathcal{H}^{m+1})^3$. The core of the proof is the convergence of the observation term  $\mathcal{J}^o(X_0^n) = \frac12 \|\xi^n(t_1)-d\|^2$ to $\frac12 \|\xi^*(t_1)-d\|^2$ where $\xi^n$ et $\xi^*$ are the Lagrangian trajectories associated to  $X_0^n$ and $X_0^*$. Indeed, we will then get the strong convergence of  $X_0^n$ to $X_0^*$ in $\mathcal{U}^{m+1}$  and in $(\mathcal{H}^{m+1})^3$ and we will easily verify that  $X_0^*$ is a minimizer of $\mathcal{J}$.\\
Let  $X^n$ be the solution of equations (\ref{2eq:1}) and (\ref{2eq:2}) associated with $X^n_0$. We choose $t^{*}$ according to  theorem \ref{2prop:1} ($t^{*}$ depends on the norm of $X_0$ and we know that the norm of $X^n_0$ is bounded, thus we can find $t^{*}$ satisfying for every $X^n$). This proposition states that the sequence  $(X^n)$ is bounded (uniformly in $n$) in $\mathcal{C}([0,t^{*}],(\mathcal{H}^{m+1})^3)$ and that the sequence $\partial_t X^n$ is uniformly bounded in $L^2(0,t^{*};(L^2_z H^m_{xy})^3)$. We then prove that  $(X^n)$ is uniformly bounded in $\mathcal{C}^{1/2}([0,t^{*}],(L^2_z H^m_{xy})^3)$:
$$ %
%\label{2eq:59}
\begin{array}{rcl}
\|X^n(t_2)-X^n(t_1)\|_{2,m} &=&  \|\int_{t_1}^{t_2} \partial_t X^n \|_{2,m} \smallbreak\\
&\leq& \sqrt{t_2-t_1} (\int_{t_1}^{t_2} \|\partial_t X^n \|^2_{2,m})^{1/2} \smallbreak\\
&\leq& \sqrt{t_2-t_1}\, C
\end{array}
$$ %
where $C$ does not depend on $n$ or any $t_i$. Thus we get:
$$ %
%\label{2eq:94}
\begin{array}{rcl}
(X^n) &\textrm{is uniformly bounded in}&\mathcal{C}([0,t^{*}],(\mathcal{H}^{m+1})^3)\\
(X^n) &\textrm{is uniformly bounded in}&\mathcal{C}^{1/2}([0,t^{*}],(L^2_z H^m_{xy})^3)\\
\Rightarrow\quad (X^n) &\textrm{is uniformly bounded in} &\mathcal{C}^{\theta/2}([0,t^{*}],([\mathcal{H}^{m+1},L^2_z H^m_{xy}]_\theta)^3)
\end{array}
$$ %
for all $\theta\in[0,1]$.\\
We now prove that there exists a small  $\theta>0$ such that, for $\delta>0$ small enough, the space $[\mathcal{H}^{m+1},L^2_z H^m_{xy}]_\theta$ is compact in the following  space $\mathcal{L}^\delta$:
$$ %
%\label{2eq:98}
\mathcal{L}^\delta = \big\{ u\in H^{1/2+\delta}_z H^{2+\delta}_{xy},\, \textrm{ periodic in } x,y,\,u=0\textrm{ on }\{z=0,z=a\}\times\mathbb{T}^2 \big\}
$$ %
For $\delta>0$, $\mathcal{L}^\delta$ is a subset of the space of the functions continuous in $z$ and Lipschitz in $(x,y)$. \\
We now describe the Hilbert interpolation  $[\mathcal{H}^{m+1},L^2_z H^m_{xy}]_\theta$, thanks to Fourier series for  $(x,y)\in\mathbb{T}^2$ and sine series in $z\in(0,a)$. Let $\zeta=(\xi,\eta)\in\mathbb{R}^2$ be the Fourier variable associated with the Fourier series on the torus, and  $\kappa\in\mathbb{N}^*$ the Fourier variable associated with sine series on  $(0,a)$.\\
We then have:
$$ %
%\label{2eq:95}
\begin{array}{rcl}
\varphi(x,y,z)\in L^2_z H^m_{xy} &\Leftrightarrow&  W_1(\zeta,\kappa)\,\hat{\varphi}(\zeta,\kappa) \in\ell^2(\mathbb{Z}^2\times\mathbb{N}^*)\\
\varphi(x,y,z)\in\mathcal{H}^{m+1} &\Leftrightarrow& W_2(\zeta,\kappa)\,\hat{\varphi}(\zeta,\kappa) \in\ell^2(\mathbb{Z}^2\times\mathbb{N}^*)
\end{array}
$$ %
where $W_1$ and $W_2$ are weights defined as follows:
$$ %
%\label{2eq:96}
\begin{array}{rcl}
  W_1(\zeta,\kappa)&=&(1+|\zeta|^2)^{\frac m2} \\
  W_2(\zeta,\kappa)&=&(1+|\zeta|^2)^{\frac{m+1}{2}}+(1+|\zeta|^2)^{\frac m2}(1+|\kappa|^2)^{\frac 12}
\end{array}
$$ %
For all $\theta>0$ small enough we get:
$$ %
%\label{2eq:97}
\varphi\in [\mathcal{H}^{m+1},L^2_z H^m_{xy}]_\theta \quad\Leftrightarrow\quad W_1^\theta W_2^{1-\theta} \hat{\varphi} \in\ell^2(\mathbb{Z}^2\times\mathbb{N}^*)
$$ %
with: 
$$ %
%\label{2eq:103}
\begin{array}{rcl}
W_1^\theta \,W_2^{1-\theta} &=& (1+|\zeta|^2)^{\frac {m\theta}2}((1+|\zeta|^2)^{\frac{m+1}{2}}+(1+|\zeta|^2)^{\frac m2}(1+|\kappa|^2)^{\frac 12})^{(1-\theta)}\\
&=& (1+|\zeta|^2)^{\frac {m}2}((1+|\zeta|^2)^{\frac{1}{2}}+(1+|\kappa|^2)^{\frac 12})^{(1-\theta)}
\end{array}
$$ %
For $\delta>0$ small enough, we also have:
$$ %
%\label{2eq:99}
\varphi\in \mathcal{L}^\delta \quad\Leftrightarrow\quad W_\delta (\zeta,\kappa)\,\hat{\varphi}(\zeta,\kappa) \in\ell^2(\mathbb{Z}^2\times\mathbb{N}^*)
$$ %
where: 
$$ %
%\label{2eq:100}
W_\delta(\zeta,\kappa)\quad=\quad(1+|\zeta|^2)^{\frac {2+\delta}2}(1+|\kappa|^2)^{\frac {\delta+1/2}2}
$$ %
To establish that the injection  $[\mathcal{H}^{m+1},L^2_z H^m_{xy}]_\theta$ is compact (for  $\theta>0$ small enough) in $\mathcal{L}^\delta$ (for $\delta>0$ small enough)  it is sufficient to find  $\delta$ and $\theta$ such that:
$$ %
%\label{2eq:101}
\lim_{|\zeta|,|\kappa|\rightarrow+\infty}\frac{W_\delta(\zeta,\kappa)}{W_1(\zeta,\kappa)^\theta \,W_2(\zeta,\kappa)^{1-\theta}} = 0
$$ %
We have:
$$ %
%\label{2eq:102}
\begin{array}{rcl}
\dfrac{W_\delta}{W_1^\theta \,W_2^{1-\theta}}&=&\dfrac{(1+|\zeta|^2)^{\frac {2+\delta}2}(1+|\kappa|^2)^{\frac {\delta+1/2}2}}   {(1+|\zeta|^2)^{\frac {m}2}((1+|\zeta|^2)^{\frac{1}{2}}+(1+|\kappa|^2)^{\frac 12})^{(1-\theta)}}\smallbreak\\
&\leq&\dfrac{(1+|\zeta|^2)^{\frac {\delta}2}(1+|\kappa|^2)^{\frac {\delta+1/2}2}}   {((1+|\zeta|^2)^{\frac{1}{2}}+(1+|\kappa|^2)^{\frac 12})^{(1-\theta)}}\smallbreak\\
&\leq& \dfrac{\frac13(1+|\zeta|^2)^{\frac{3\delta}2} + \frac23(1+|\kappa|^2)^{\frac {3\delta+3/2}4}}   {((1+|\zeta|^2)^{\frac{1}{2}}+(1+|\kappa|^2)^{\frac 12})^{(1-\theta)}}
\end{array}
$$ %
using that $m\geq 2$ and $ab\leq \frac{a^p}p+\frac{b^q}q$ for $a,b$ positive and $p=3$, $q=\frac32$. Then we get:
$$ %
\label{2eq:104}
\begin{array}{rcl}
\dfrac{W_\delta}{W_1^\theta \,W_2^{1-\theta}}&\leq&(1+|\zeta|^2)^{\frac {3\delta}2-\frac{1-\theta}2}    + (1+|\kappa|^2)^{\frac {3\delta+3/2}4-\frac{1-\theta}2}\smallbreak\\
\end{array}
$$ %
For $\theta$ and $\delta$ positive and small enough we have:
$$ %
%\label{2eq:105}
\frac {3\delta}2-\frac{1-\theta}2 \simeq -\frac12 \quad\textrm{ et }\quad \frac {3\delta+3/2}4-\frac{1-\theta}2\simeq -\frac18
$$ %
Then we can find $\theta$ and $\delta$  positive and small enough such that the right hand side of equation (\ref{2eq:104}) converges to 0. We have finally shown that  $(X^n)$ is uniformly bounded in $\mathcal{C}^{\theta/2}([0,t^{*}],([\mathcal{H}^{m+1},L^2_z H^m_{xy}]_\theta)^3)$ with  $[\mathcal{H}^{m+1},L^2_z H^m_{xy}]_\theta$ compact in $\mathcal{L}^\delta$. 
% ici on utilise le resultat suivant :
% si Y est compact dans X, alors pour tout eps on peut ecrire X comme somme directe de E_esp et F_eps
% avec E de dim finie et pour tout y ds Y, si on decompose y=y_E+y_F
% alors la norme de y_F dans X est inf a eps fois la norme de y dans Y
% grace a ca on peut montrer que la suite X^n converge bien la ou il faut.
Then $(X^n)$ converges strongly to  $X^*$ in $\mathcal{C}([0,t^{*}],(\mathcal{L}^\delta)^3)$. And finally we have:
$$ %
%\label{2eq:106}
\mathcal{J}^o(X_0^n) = \frac12 \|\xi^n(t_1)-d\|^2 \quad\rightarrow\quad \frac12 \|\xi^*(t_1)-d\|^2 = \mathcal{J}^o(X_0^*)
$$ %
Then we prove that $(X_0^n)$ converges strongly to $X_0^*$ in $\mathcal{U}^{m+1}$:
$$ %
%\label{2eq:107}
\begin{array}{rcl}
\mathcal{J}(X_0^*)&=&\|X_0^*\|^2_{\mathcal{U}^{m+1}}+\mathcal{J}^o(X_0^*)\\
&\leq& \underline{\lim}\|X_0^n\|^2_{\mathcal{U}^{m+1}}+\lim\mathcal{J}^o(X_0^n)\\
&\leq& \inf_{X_0\in\mathcal{U}^{m+1}} \mathcal{J}(X_0) 
\end{array}
$$ %
thus $\mathcal{J}(X_0^*) = \inf_{X_0\in\mathcal{U}^{m+1}} \mathcal{J}(X_0)$, then $\|X_0^n\|_{\mathcal{U}^{m+1}} \rightarrow \|X_0^*\|_{\mathcal{U}^{m+1}}$, therefore $(X_0^n)$ converges strongly to  $X_0^*$ in $\mathcal{U}^{m+1}$, and $X_0^*$ is a minimizer of $\mathcal{J}$.\\
The strong convergence of $(X^n)$ in $\mathcal{C}([0,t^{*}],(\mathcal{L}^\delta)^3)$ implies also:
$$ %
%\label{2eq:108}
X^*(t=0)=X_0^*
$$ %

\subsection{Passage to the limit in the state equation}
In this section we prove that the limit  $X^*$ satisfies equation (\ref{2eq:1}) and boundary conditions (\ref{2eq:2}). To do so, we first pass to the limit in the weak formulation of the equation. Let $X' = (u', v', \theta') \in (\mathcal{D}((0,t^{*})\times\Omega))^3$ be a test function satisfying conditions (\ref{2eq:2}). We now pass to the limit in the five following terms:
$$ %
\begin{array}{rcl}
%\label{2eq:60}
T_1 &=& \iint \partial_t u^n u' + \partial_t v^n v' + K \partial_t \theta^n \theta'\\
T_2 &=& \iint -\nu \Delta u^n u'  -\nu \Delta v^n v'  -\nu K \Delta \theta^n \theta'\\
T_3 &=& \iint -\alpha v^n u' + \alpha u^n v'\\
T_4 &=& \iint -\beta \theta^n w' + K \gamma w^n \theta'\\
T_5 &=& \iint (U^n.\nabla_2 u^n +w^n \partial_z u^n) u' + (U^n.\nabla_2 v^n +w^n \partial_z v^n) v'\\
&&\qquad  + K (U^n.\nabla_2 \theta^n +w^n \partial_z \theta^n) \theta'
\end{array}
$$ %
$T_1$. As we proved that $(\partial_t X^n)$ is uniformly bounded in $L^2(0,t^{*};(L^2_z H^m_{xy})^3)$, it converges thus weakly in $L^2(0,t^{*};(L^2_z H^m_{xy})^3)$. We can prove easily that its limit is  $\partial_t X^*$, e. g., for $\partial_t u^n$:
$$ %
%\label{2eq:61}
\iint \partial_t u^n u' = - \iint u^n \partial_t u' \,\, \rightarrow \,\, - \iint u^* \partial_t u' = \iint \partial_t u^* u'
$$ %
Then we get for $T_1$:
$$ %
%\label{2eq:62}
\begin{array}{rl}
T_1 \rightarrow &\iint \partial_t u^* u' + \partial_t v^* v' + K \partial_t \theta^* \theta'
\end{array}
$$ %
$T_2$. A space integration by parts gives:
$$ %
%\label{2eq:63}
\iint -\nu \Delta u^n u' = -\iint \nu u^n \Delta u'
$$ %
To the limit we obtain:
$$ %
%\label{2eq:64}
\begin{array}{rl}
T_2 \rightarrow &-\iint \nu u^* \Delta u' -\iint \nu v^* \Delta v' -K \iint \nu \theta^* \Delta \theta' 
\end{array}
$$ %
$T_3$. We easily have:
$$ %
%\label{2eq:65}
\begin{array}{rl}
T_3 \rightarrow &\iint -\alpha v^* u' + \alpha u^* v'
\end{array}
$$ %
$T_4$. The first part is easy. For the second we have:
$$ %
%\label{2eq:66}
\begin{array}{rcl}
\iint   w^n \theta' &=& -\iint \big[ \int_0^z (\partial_x u^n(z') + \partial_y v^n(z'))\,dz' \big] \theta' \smallbreak\\
&=& \int_0^{t^{*}} \int_{xy} \int_{z'=0}^a  (\partial_x u^n(z') + \partial_y v^n(z'))\big[ \int_{z=z'}^a  \theta'(z) \,dz \big] \,dz'\smallbreak\\
&=& \iint (\partial_x u^n + \partial_y v^n)\big[ \int_{z'=z}^a  \theta'(z') \,dz' \big]
\end{array}
$$ %
thanks to Fubini's equality. As $(X^n)$ converges strongly in $L^2(0,t^{*};(L^2_z H^m_{xy})^3)$ with $m\geq 2$, we can pass to the limit, using Fubini again and posing $w^*=\int_0^z (\partial_x u^*(z') + \partial_y v^*(z'))\,dz'$:
$$ %
%\label{2eq:67}
\begin{array}{rl}
T_4 \rightarrow &-\beta\iint  \theta^* w' - K \gamma \iint w^* \theta'
\end{array}
$$ %
$T_5$. Let us consider the following part (the other parts can be evaluated similarly):
$$ %
%\label{2eq:68}
\begin{array}{ccccc}
\iint (w^n\partial_z u^n - w^*\partial_z u^*) u'& =& \iint w^*(\partial_z u^n -\partial_z u^*)u'& +& \iint (w^n-w^*)\partial_z u^n u'\smallbreak\\
&=& T_{5,1} &+&T_{5,2}
\end{array}
$$ %
We know that $\partial_z u^n -\partial_z u^*$ converges weakly to $0$ in $L^2(0,t^{*};L^2_{z}L^2_{xy})$. We also have  $w^* \in L^2(0,t^{*};L^2_{z}H^2_{xy})$ because $X^* \in L^2(0,t^{*};\mathcal{U}^{m+1})$. We now prove that $w^*u' \in L^2(0,t^{*};L^2_{z}L^2_{xy})$ as in the proof of proposition \ref{2lem:2}:
$$ %
%\label{2eq:69}
\begin{array}{rcl}
\iint  |w^* u'|^2 &\leq& \int_{t} \int_z \|w^* u'\|_{H^2_{xy}}^2\smallbreak\\
&\leq&  \int_{t} \int_z \|w^*\|_{H^2_{xy}}^2 \|u'\|_{H^2_{xy}}^2\smallbreak\\
&\leq&  \int_{t} \|u'\|_{L^2_zH^2_{xy}} \|\partial_z u'\|_{L^2_zH^2_{xy}} \|w^*\|_{L^2_zH^2_{xy}}^2\smallbreak\\
&\leq& \sup_t \big[ \, \|u'\|_{L^2_zH^2_{xy}} \|\partial_z u'\|_{L^2_zH^2_{xy}} \, \big] \, \|w^*\|_{L^2_tL^2_zH^2_{xy}}^2\smallbreak\\
&<& \infty
\end{array}
$$ %
Thus $T_{5,1}$ converges to 0.\\
For $T_{5,2}$ we have:
$$ %
%\label{2eq:70}
\begin{array}{rcl}
|\iint (w^n-w^*)\partial_z u^n u'|^2 &\leq& \|w^n-w^*\|^2_{L^2((0,t^{*})\times\Omega)} \|\partial_z u^n u'\|^2_{L^2((0,t^{*})\times\Omega)}\smallbreak\\
&\leq& C \|X^n-X^*\|^2_{L^2(0,t^{*};L^2_zH^1_{xy})} \int_t \int_z \|\partial_z u^n\|_{H^2_{xy}}^2 \|u'\|_{H^2_{xy}}^2\smallbreak\\
&\leq& C \|X^n-X^*\|^2_{L^2(0,t^{*};L^2_zH^1_{xy})} \|\partial_z u^n\|_{L^2(0,t^{*};L^2_zH^2_{xy})}^2\smallbreak\\
&\leq& C \|X^n-X^*\|^2_{L^2(0,t^{*};L^2_zH^1_{xy})} \|X^n\|_{L^2(0,t^{*};\mathcal{U}^{m+1})}^2\smallbreak\\
&\leq& C \|X^n-X^*\|^2_{L^2(0,t^{*};L^2_zH^1_{xy})}
\end{array}
$$ %
so $T_{5,2}$ converges also to 0.
This gives the expected limit for $T_5$:
$$ %
%\label{2eq:71}
\begin{array}{rl}
T_5 \rightarrow &\iint (U^*.\nabla_2 u^* +w^* \partial_z u^*) u' + (U^*.\nabla_2 v^* +w^* \partial_z v^*) v'\\
&\qquad  + K (U^*.\nabla_2 \theta^* +w^* \partial_z \theta^*) \theta'\\
\end{array}
$$ %
Then the limit $X^*$ satisfies the following equation, for all $X' \in (\mathcal{D}((0,t^{*})\times\Omega))^3$ satisfying (\ref{2eq:2}):
\begin{equation}
\label{2eq:72}
\begin{array}{rl}
0=&  \iint \big(\partial_t u^*  u' +\partial_t v^*   v'+\partial_t \theta^*  \theta' \big) \smallbreak\\
&-\iint (\nu u^* \Delta u' + \nu v^* \Delta v' +K \nu \theta^* \Delta \theta') \smallbreak\\
&+\iint (-\alpha v^* u' + \alpha u^* v')\smallbreak\\
&+\iint( -\beta \theta^* w' - K \gamma \big[ \int_0^z (\partial_x u^*(z') + \partial_y v^*(z'))\,dz' \big] \theta' )\smallbreak\\
&- \iint \big(\,\, (U^*.\nabla_2 u^* +w^* \partial_z u^*) u' + (U^*.\nabla_2 v^* +w^* \partial_z v^*) v'\\
&\qquad  + K (U^*.\nabla_2 \theta^* +w^* \partial_z \theta^*) \theta' \,\,\big)
\end{array}
\end{equation}
We now prove that this defines a linear functional which is continuous on $L^2(0,t^{*};(H^1_0(\Omega))^3$. It is obvious that the terms involving linearly $X^*$ and its derivates define a linear functional. It remains to study terms involving  $X^*$ non linearly. Let us define:
$$ %
%\label{2eq:109}
\phi(u') = \iint  (u^* \partial_x u^* +w^* \partial_z u^*) u'
$$ %
We will prove that $\phi$ is continuous on $L^2(0,t^{*};L^2(\Omega))$ and therefore on $L^2(0,t^{*};H^1_0(\Omega))$. As we have $X^* \in \mathcal{C}(0,t^{*},\mathcal{U}^{m+1})$, every space derivative of  $u^*$, and also $w^*$ and $\partial_z w^*$, are in $\mathcal{C}(0,t^{*},L^2_z H^2_{xy})$. Then we have, as in the proof of proposition \ref{2lem:2}:
$$ %
%\label{2eq:110}
\begin{array}{rcl}
\iint  |u^* \partial_x u^* +w^* \partial_z u^*|^2 &=& \int_{t,z}  \|u^* \partial_x u^* +w^* \partial_z u^*\|^2_{L^2_{xy}}\\
&\leq& C\int_{t,z}  \|u^* \partial_x u^* +w^* \partial_z u^*\|^2_{H^2_{xy}}\\
&\leq& C\int_{t,z}  \|u^*\|^2_{H^2_{xy}}\|\partial_x u^*\|^2_{H^2_{xy}} +\|w^* \|^2_{H^2_{xy}}\|\partial_z u^*\|^2_{H^2_{xy}}\\
&\leq&C \int_{t} \|u^*\|_{L^2_z H^2_{xy}} \|\partial_z u^*\|_{L^2_z H^2_{xy}}\|\partial_x u^*\|^2_{H^2_{xy}} \\
&&  + C \int_{t} \|w^* \|_{L^2 H^2_{xy}}  \|\partial_z w^* \|_{L^2 H^2_{xy}}    \|\partial_z u^*\|^2_{H^2_{xy}}\\
&\leq& C \sup_t \|X^*(t)\|^4_{\mathcal{U}^{m+1}}
\end{array}
$$ %
Thus $\phi$ is a continuous linear functional on $L^2(0,t^{*};H^1_0(\Omega))$. Similarly we deal with every other non linear term in $X^*$.\\
Finally, formula (\ref{2eq:72}) defines a linear functional, which is continuous on  $L^2(0,t^{*};H^1_0(\Omega))^3$  and vanishes on the elements of $L^2(0,t^{*};(H^1_0(\Omega))^3$ satisfying the boundary conditions (\ref{2eq:2}). As in  \cite{TemamZiane04} and in the theory of Navier-Stokes equations (e.g., see \cite{Lions69,Temam77}), we can reintroduce the pressure: there exists $p^* \in \mathcal{D}'(0,t^{*},L^2(\Omega))$ such that $(X^*,p^*)$ is a solution of the Primitive Equations  (\ref{2eq:1}) with boundary conditions (\ref{2eq:2}).

% 
%%-----------------------------
\section{Numerical experiments}
\label{sec:num}
%%-----------------------------
%
In this section, we briefly present some illustrative numerical results. A detailed description of the numerical setup and more in-depth results can be found in \cite{Nodet06}.%
\subsection{Numerical setup}
This problem was addressed with a realistic state-of-the-art Primitive Equations ocean model, namely OPA code, developed by LODYC (see \cite{Madec99}). The model is set-up in a classical double-gyre wind-driven configuration: it is representative of mid-latitude ocean circulation, where a non-linear and non-stationary jet-stream (such as the Gulf Stream) develops at the convergence of the subpolar gyre and the subtropical gyre.\\

The control optimal problem that we solved numerically is very similar to the one presented before, the regularization term in the cost function (\ref{2eq:6}) being different:
\begin{equation}
\begin{array}{ccccc}
\label{numeq:6}
\mathcal{J}(X_0) & = & \frac{1}{2} \sum_{j=1}^M \sum_{i=1}^N\|\xi_j (t_i) -d_i^j\|^2 &+ & \frac{\omega}{2} \| X_0 -X_b\|_{B}^2\\
& =& \mathcal{J}^o(X_0)  &+&\omega \, \mathcal{J}^b(X_0)  
\end{array} 
\end{equation}
The observation term $\mathcal{J}^o$ is the same, except we have $M$ floats drifting and $N$ time-sampling of their positions. The background term $\mathcal{J}^b$ involved a so-called background state $X_b$, containing a priori information (such as climatology, or results of a previous forecast) and a background error covariance matrix $B$ to define the $B$-norm:
$$ %
\|X\|_B^2 = \, ^tX B^{-1} X
$$ %

\begin{figure}
\includegraphics[angle=-90,width=\textwidth]{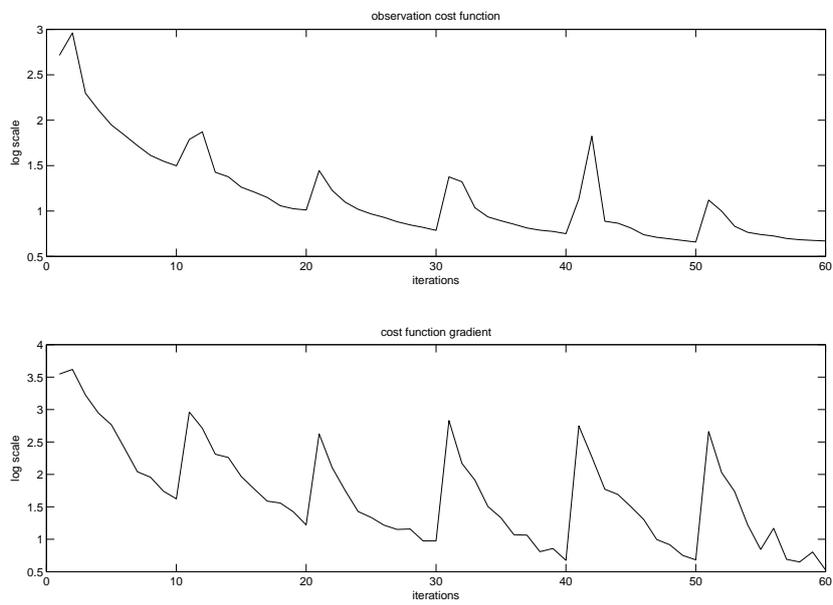}
\caption{Observation term of the cost function (top) and gradient of the total cost function (bottom) as a function of the iterations number. The linearization re-occurs every ten iterations.
\label{fig:1} 
}
\end{figure}

\begin{figure}
\includegraphics[angle=-90,width=\textwidth]{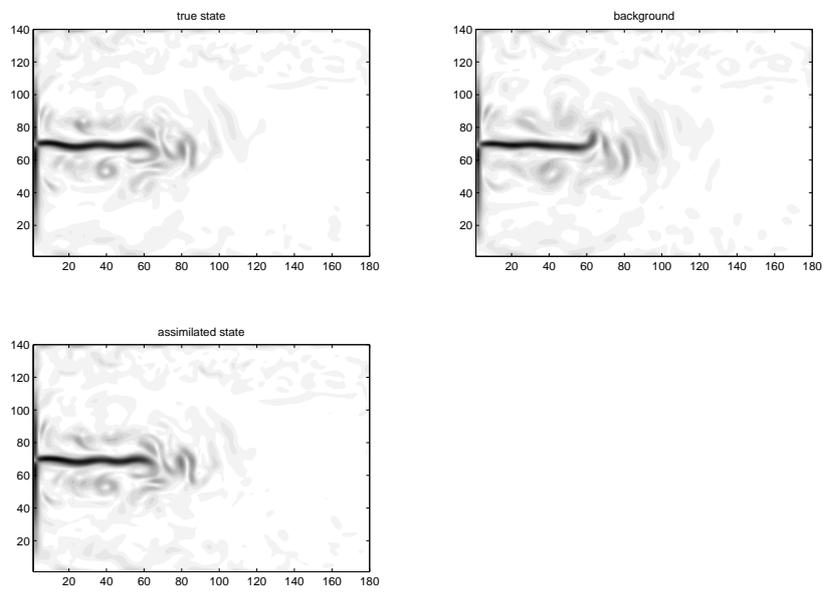}
\caption{Kinetic energy of the true state, background (no assimilation) and assimilated state at the surface after thirty days; the x-coordinate represent longitude grid-points and y-coordinate are latitude grid-points. 
\label{fig:2} 
}
\end{figure}

\begin{figure}
\includegraphics[angle=-90,width=\textwidth]{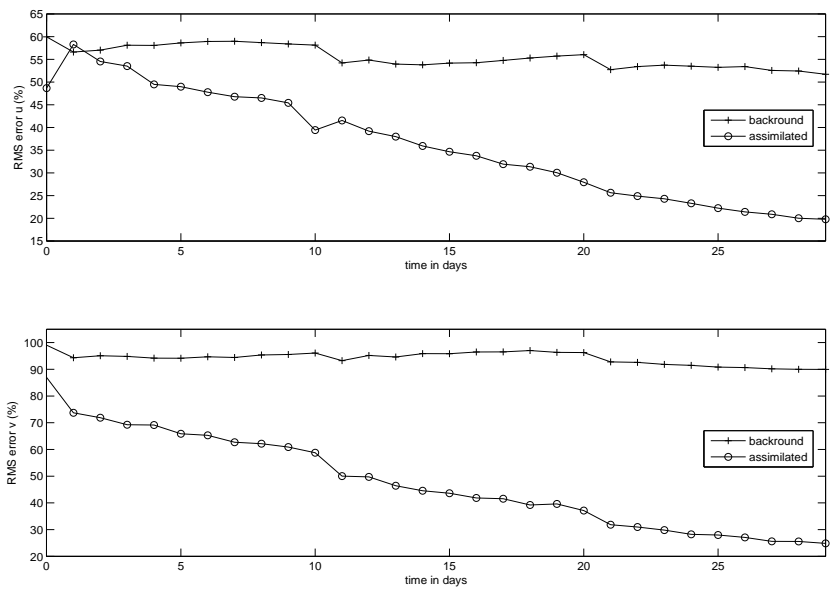}
\caption{Relative RMS errors for $u$ (top) and $v$ (bottom) on the whole grid as a function of time for the background and the assimilated state, as a function of the time.
\label{fig:3} 
}
\end{figure}

\subsection{Four-dimensional variational data assimilation}
The cost function (\ref{numeq:6}) is minimized thanks to an iterative process involving its gradient. Said gradient is computed thanks to the adjoint of the ocean model, and also the adjoint of the observation operator (\ref{2eq:121}). This method, called four-dimensional variational data assimilation, has been introduced in meteorology by Le Dimet (see \cite{LeDimet82,LeDimet86}) and is based on \cite{Lions68}. Although this method has been introduced for linear models and observation operators, it can be extended to ``reasonably" non-linear cases: in \cite{Courtier94} the authors proposed an incremental approach, which consists in linearizing the operators around a given trajectory and then proceed to minimize the quadratic cost function. Once the minimum is reached, the operators are once again linearized around it, and the new cost function is minimized, and so on. Figure \ref{fig:1} presents the decrease of the cost function and its gradient as a function of the iteration number. Every ten iterations, the linearized operators and their adjoints are updated.

%fig1
%
\subsection{Illustration}
The test-experiment presented here involves $M=1000$ floats, drifting at 1000 meters depth, whose positions are sampled every day during ten days. In order to achieve convergence, we had to assume that we have a priori information on temperature and salinity, so that the assimilation process aims to reconstruct the velocities, from the positions information. We performed identical twin experiments: a given output of the ocean model is called ``true state" and is used to generate observations. Then the iterative process is initialized with the background (which is equal to the true state except it has wrong velocities), and the observations are assimilated. We then get an ``assimilated stated" which should be close to the true state. The experiment presented here consists of three assimilation process on three successive ten-days periods.\\
Figure \ref{fig:2} presents the square root of the kinetic energy at time $T=30$ days, at the surface, for the true state (reference), the background (no assimilation) and the assimilated state. The assimilation process reconstructed very well the true state. Let us note that the horizontal velocity field presented in Figure \ref{fig:2} is the surface one, whereas the floats drift at 1000 meters depth. The assimilation thus transfered information to every vertical level and not only the 1000 meters deep one.\\
%
%fig2
%
Figure \ref{fig:3} presents the evolution of the relative RMS error (with respect to the true state) over the 30 days time window. For any velocity field $(u(x,y,z,t),v(x,y,z,t))$, such relative RMS error for $u$ is a function of time, computed as follows:
$$ %
%\label{eq:err}
\mathcal{E}(u;t) = \Big( \frac{\int_\Omega |u_t(x,y,z,t) - u(x,y,z,t)|^2 \,dx\,dy\,dz }{\int_\Omega |u_t(x,y,z,t)|^2 \,dx\,dy\,dz}\Big)^{1/2}
$$ %
where $u_{t}$ is the true velocity. Similarly we can compute $\mathcal{E}(v;t)$. We then can compare the errors for the background $\mathcal{E}(u_{b};t)$ and the assimilated state $\mathcal{E}(u_{a};t)$. We can see that after thirty days the error has been divided by 3 for $u$ and $4$ for $v$.

%fig3

%%-----------------------------
%%      your bibliography
\bibliographystyle{plain}
\bibliography{biblio_cocv}
%%-----------------------------

\end{document}